\documentclass[11pt]{amsart}                                                                                                                                           
                                                                                                                                                                       
\usepackage{amsmath}                                                                                                                                                   
\usepackage{amssymb}                                                                                                                                                   
\usepackage{latexsym}                                                                                                                                                  
\usepackage{amsfonts}    
\begin{document}                                                                                                                                                       
%Commands Used%                                                                                                                                                        
                                                                                                                                                                       
\newcommand{\ci}[1]{_{ {}_{\scriptstyle #1}}}

\newcommand{\norm}[1]{\ensuremath{\|#1\|}}                                                                                                                             
\newcommand{\abs}[1]{\ensuremath{\vert#1\vert}}                                                                                                                        
\newcommand{\p}{\ensuremath{\partial}}

\newcommand{\pbar}{\ensuremath{\bar{\partial}}}                                                                                                                        
\newcommand{\db}{\overline\partial}                                                                                                                                    
\newcommand{\D}{\mathbb{D}}                                                                                                                                            
\newcommand{\T}{\mathbb{T}}                                                                                                                                            
\newcommand{\C}{\mathbb{C}}                                                                                                                                            
\newcommand{\N}{\mathbb{N}}                                                                                                                                            
\newcommand{\td}{\widetilde\Delta}

\newcommand{\La}{\langle }                                                                                                                                             
\newcommand{\Ra}{\rangle }                                                                                                                                             
\newcommand{\tr}{\operatorname{tr}}                                                                                                                                    
\newcommand{\ran}{\operatorname{Ran}}                                                                                                                                  
\newcommand{\vf}{\varphi}                                                                                                                                              
\newcommand{\e}{\varepsilon}                                                                                                                                           
                                                                                                                                                                       
\newcommand{\f}[2]{\ensuremath{\frac{#1}{#2}}}

\newcommand{\be}{\mathbf{e}}                                                                                                                                           
\newcommand{\clos}{\operatorname{clos}}                                                                                                                                
\newcommand{\rank}{\operatorname{rank}}                                                                                                                                
                                                                                                                                                                       
\newcommand{\bz}{\mathbf{z}}

\newcommand{\tto}{\!\!\to\!}                                                                                                                                           
\newcommand{\wt}{\widetilde}                                                                                                                                           
\newcommand{\shto}{\raisebox{.3ex}{$\scriptscriptstyle\rightarrow$}\!}                                                                                                 
\newcommand{\bL}{\mathbf{L}}                                                                                                                                                                      
%%%%%%%%%%%%%%%%%%%%%%%%%%%%                                                                                                                                           
                                                                                                                                                                       
\newcommand{\entrylabel}[1]{\mbox{#1}\hfill}                                                                                                                           
                                                                                                                 %%%%%%%%%%%%%%%%%%%%

\makeatletter
 
\def\multilimits@{\bgroup
  \Let@
  \restore@math@cr
  \default@tag
 \baselineskip\fontdimen10 \scriptfont\tw@
 \advance\baselineskip\fontdimen12 \scriptfont\tw@
 \lineskip\thr@@\fontdimen8 \scriptfont\thr@@
 \lineskiplimit\lineskip
 \vbox\bgroup\ialign\bgroup\hfil$\m@th\scriptstyle{##}$\hfil\crcr}
\def\Sb{_\multilimits@}
\def\Sp{^\multilimits@}
\def\endSb{\crcr\egroup\egroup\egroup}
\let\endSp=\endSb

\makeatother
%%%%%%%%%%%%%%%%%%%%%%%%%%%%%%%%                                                     
                                                                                                                                                                       
\newenvironment{entry}
{\begin{list}{X}%
  {\renewcommand{\makelabel}{\entrylabel}%
      \setlength{\labelwidth}{55pt}%
      \setlength{\leftmargin}{\labelwidth}%\labelsep}%
      \addtolength{\leftmargin}{\labelsep}%
   }%
}%
{\end{list}}

%%%%%%%%%%%%%%%%%%%%%%%%%%%%                                                                                                                                           
                                                                                                                                                                       
%\newenvironment{theorem}{\newline\begin{sloppypar}\noindent\textbf{Theorem}}{\hspace*{\fill}\                                                                         
%\end{sloppypar}}                                                                                                                                                      
                                                                                                                                                                       
\numberwithin{equation}{section}                                                                                                                                       
                                                                                                                                                                       
\newtheorem{thm}{Theorem}[section]                                                                                                                                     
\newtheorem{lm}[thm]{Lemma}                                                                                                                                            
\newtheorem{cor}[thm]{Corollary}                                                                                                                                       
\newtheorem{prop}[thm]{Proposition}

\theoremstyle{remark}                                                                                                                                                  
\newtheorem{rem}[thm]{Remark}                                                                                                                                          
\newtheorem*{rem*}{Remark}

\title{The Matrix-Valued $H^{p}$ Corona Problem in the Disk and                                                                                                        
Polydisk}                                                                                                                                                              
\author{Sergei Treil \and Brett D. Wick}                                                                                                                               
                                                                                                                                                                       
\begin{abstract}                                                                                                                                                       
In this paper we consider the matrix-valued $H^p$ corona problem in                                           
the disk and polydisk.  The result for the disk is rather well known,                                                                                                  and is usually obtained from the classical Carleson Corona Theorem by                                                                                                  
linear algebra.  Our proof provides a streamlined way of obtaining                                                                                                     this result and allows one to get a better estimate on the norm of                                                                                                     
the solution.  In particular, we were able to improve the estimate                                                                                                     
found in the recent work of T. Trent in \cite{Trent}.  Note that, the                                                                                                  
solution of the $H^{\infty}$ matrix corona problem in the disk can be                                                                                                  
easily obtained from the $H^{2}$ corona problem either by                                                                                                              
factorization, or by the Commutant Lifting Theorem.                                                                                                                    
The $H^{p}$ corona problem in the polydisk was originally solved by                                                                                                    
K.C. Lin in \cite{Lin1} and \cite{Lin2}.  The solution used Koszul                                                                                                     
complexes and was rather complicated because one had to consider                                                                                                       
higher order $\pbar$-equations.  Our proof is more transparent and it                                                                                                  
improves upon Lin's result in several ways.  First we were able to                                                                                                     
show that the norm of the solution is independent of the number of                                                                                                     
generators.  Additionally, we illustrate that the norm of the                                                                                                         
solution grows at most proportionally to the dimension of the polydisk.  Our approach is based on one that was                                                                                                      
originated by M. Andersson in \cite{Andersson1}.  In the disk it                                                                                                       
essentially depends on Green's Theorem and duality to obtain the                                                                                                       
estimate.  In the polydisk we use Riesz projections to reduce the                                                                                                   
problem to the disk case.                                                                                                                                              
\end{abstract}                                                                                                                                                         
                                                                                                                                                                       
\maketitle                                                                                                                                                             
\setcounter{tocdepth}{1}                                                                                                                                                  
\tableofcontents                                                                                                                                                       
\section*{Notation}                 

%\addcontentsline{toc}{section}{Notation}                                                                                                                              
\begin{entry}                                                                                                                                                          
\item[$:=$] equal by definition;\\                                                                                                                                     
\item[$\C$] the complex plane;\\                                                                                                                                       
\item[$\D$] the unit disk, $\D:=\{z\in\C:\abs{z}<1\}$;\\                                                                                                               
\item[$\T$] the unit circle, $\T:=\p\D=\{z\in\C:\abs{z}=1\}$;\\                                                                                                        
\item[$F^{*}$] the conjugate transpose of the matrix $F$;\\                                                                                                            
\item[$d\mu$] measure on $\D$ with                                                                               $d\mu=\f{2}{\pi}\log\f{1}{\abs{z}}dxdy$;\\                                                                                                                             
\item[$dm$] normalized Lebesgue measure on $\T$;\\                                                                                                                     
\item[$\La\cdot ,\cdot\Ra$] standard inner product in $\C^{n}$;\\                                                                                                      
\item[$\norm{\cdot}$]  norm; since we are dealing with matrix and                                                                                                      
operator-valued functions this symbol is a bit overloaded, but we hope it                                                                                                      
will not cause any confusion. The norm in the function spaces can be                                                                                                   
always distinguished by subscript. Thus for a vector-valued function                                                                                                   
$f$ the symbol $\|f\|_2$ denotes its $L^2$-norm, but the symbol                                                                                                        
$\|f\|$ stands for the scalar valued function whose value at a point                                                                                                   
$z$ is the norm of the vector $f(z)$; \\                                                                                                                               
                                                                                                                                                                       
\item[$\|f\|_2$] $L^2$-norm of the vector-valued function $f$ on the                                                                                                   
unit circle, \newline $\displaystyle \|f\|_2^2 := \int_{\T}                                                                                                            
\|f(z)\|^2 dm(z)$;                                                                                                                                                     
                                                                                                                                                                       
\item[$\norm{F}$] norm of $F$ as an operator;\\                                                                                                                        
\item[$H^{\infty}(\D)$] space of bounded analytic functions on $\D$                                                                                                    
with the supremum norm;\\                                                                                                                                              
\item[$H^{\infty}(\D;E\to E_*)$] operator Hardy class of bounded                                                                                                       
analytic functions from the disk whose values are bounded operators                                                                                                    
from $E$ to $E_*$, \newline                                                                                                                                            
$\displaystyle\norm{F}_{\infty}:=\sup_{z\in\D}\norm{F(z)}$;                                                                                                                                                                                                                   
\end{entry}

We will need some specialized notation when dealing with the matrix-valued Corona problem in the polydisk.  The reason for the additional notation is that  we must keep track of the extra variables that are present.  The following notation only applies to the last two sections.

\begin{entry}
\item[$\C^{n}$] n-dimensional complex space;\\
\item[$\D^{n}$] the polydisk in $\C^{n}$, 
$\D^{n}:=\{\bz\in\C^{n}:\abs{z_{i}}<1\,\forall i=1,\ldots,n\}$;\\
\item[$\T^{n}$] the \v{S}ilov boundary of $\D^{n}$, 
$\T^{n}:=\p_{S}\D^{n}=\{\bz\in\C^{n}:\abs{z_{i}}=1\,\forall 
i=1,\ldots,n\}$;\\
\item[$\p_{i}$] derivative with respect to the variable $z_{i}$;\\
\item[$\pbar_{i}$] derivative with respect to the variable 
$\bar{z_{i}}$;\\
\item[$d\mu(z_k)$] measure on $\D$ with 
$d\mu(z_k)=\f{2}{\pi}\log\f{1}{\abs{z_{k}}}dx_{k}dy_{k}$;\\
\item[$dm(z_k)$] normalized Lebesgue measure on $\T$ with respect to 
the $z_{k}$ variable;\\
\item[$dm_{n}(\bz)$] $n$-dimensional normalized Lebesgue measure on 
$\T^{n}$, i.e., $dm_{n}(\bz)=dm(z_1)\cdots dm(z_n)$;\\
\item[$dm_{n-1}(\bz_j)$] normalized $(n-1)$-dimensional Lebesgue 
measure on $\T^{n-1}$ with $dm(z_j)$ removed, i.e., 
$dm_{n-1}(\bz_j)=dm(z_1)\cdots\widehat{dm(z_j)}\cdots dm(z_n)$, with 
the circumflex denoting omission of that term;\\
\item[$\La\cdot ,\cdot\Ra$] standard inner product in $\C^{k}$;\\
\item[$(\cdot ,\cdot)$] inner product in a Hilbert space;\\
\item[$\|f\|_2$] $L^2$-norm of the vector-valued function $f$ on the 
\v{S}ilov boundary $\T^{n}$, \newline $\displaystyle \|f\|_2^2 := 
\int_{\T^{n}} \|f(\bz)\|^2 dm_{n}(\bz)$;\\
\item[$H^{\infty}(\D^{n};E\to E_*)$] operator Hardy class of bounded                                                                                                   
analytic functions from the polydisk whose values are bounded                                                                                                          
operators from $E$ to $E_*$, \newline                                                                                                                                  
$\displaystyle\norm{F}_{\infty}:=\sup_{\bz\in\D^{n}}\norm{F(\bz)}$; 
\end{entry}

Throughout the paper all Hilbert spaces are assumed to be separable. 
We
always assume that in any Hilbert space an orthonormal basis is 
fixed, so
an operator $A:E\to E_*$ can be identified with its matrix. Thus 
besides
the usual involution $A\mapsto A^*$ ($A^*$ is the adjoint of $A$), we 
have two more: $A\mapsto A^T$ (transpose of the matrix) and $A\mapsto
\overline A$ (complex conjugation of the matrix), so $A^* =(\overline
A)^T =\overline{A^T}$. Although everything in the paper can be 
presented in an
invariant, ``coordinate-free'', form, use of transposition and complex
conjugation makes the notation easier and more transparent. 

\setcounter{section}{-1}

\section{Introduction and main result}
The classical Carleson Corona Theorem, see \cite{Carleson}, states that 
if functions
$f_{j}\in H^{\infty}(\D)$ are such that
$\sum_{j=1}^{\infty}\abs{f_{j}}^{2}\geq\delta^2>0$ then  there exist
functions $g_{j}\in H^{\infty}(\D)$ such that
$\sum_{j=1}^{\infty}g_{j}f_{j}=1$.   This is equivalent to the fact that the unit disk $\D$ is dense in the
maximal ideal space of the algebra $H^\infty$, but the importance of 
the
Corona Theorem goes much beyond the theory of maximal ideals of
$H^\infty$. 

The Corona Theorem, and especially its generalization, the so called
Matrix (Operator) Corona Theorem play an important role in 
operator
theory (such as the angles between invariant subspaces, unconditionally convergent
spectral decomposition, computation of spectrum, etc.). The Matrix 
Corona
Theorem says that if  $F\in H^\infty(\D; E_*\tto E)$ is a bounded 
analytic
function whose values are operators from a Hilbert space $E_*$, $\dim
E_*<+\infty$, to another Hilbert space $E$ such that 
\begin{equation}
\label{C}
\tag{C}
F^*(z)F(z)\ge\delta^2 I >0, \qquad \forall z\in\D, 
\end{equation}
then $F$ has a bounded analytic left inverse $G\in H^\infty 
(\D;E_*\tto
E)$, $GF\equiv I$. We should emphasize that the requirement
$\dim E_*<+\infty$ is essential here.   It was shown in \cite{Tre-cor-cex}, see
also
\cite{Treil-OTAA-89} or \cite{Tre-rank-one}, that the Operator Corona Theorem 
fails
if $\dim E_*=+\infty$. Note also that the above condition \eqref{C} is
necessary for the existence of a bounded left inverse. 

The classical Carleson Corona Theorem is a particualr case of the 
matrix
one: one just needs to consider $F$ being the column $F=(f_1, f_2, 
\ldots,
f_n)^T$. It also worth noticing that the Matrix Corona Theorem follows
from the classical one.  Using a simple linear algebra argument 
P.~Fuhrmann,
see \cite{Fuhr-corona}, was able to get  the matrix version ($\dim E_*, 
\dim
E<+\infty$) of the theorem from the classical result of Carleson. 
Later,
using the ideas from T.~Wolff's proof of the Corona Theorem 
M.~Rosenblum, 
V.~Tolokonnikov and A.~Uchiyama independently extended the Corona Theorem 
to
infinitely many functions $f_k$. Using their result, V.~Vasyunin was 
able
to get the Operator Corona Theorem in the case $\dim E_*<+\infty$, 
$\dim
E=+\infty$.

Since the Corona Theorem turns out to be very important in operator
theory, there were some attempts  to prove it using operator
methods. While these attempts were not completely successful, some 
interesting relations were discovered. In particular, it was shown 
that a function $F\in H^\infty =H^\infty (\D;E_*\tto E)$ is left 
invertible in $H^\infty$ if and only if the Toeplitz operator 
$T_{\overline F}$ is left invertible; 
here $\overline F$ denotes the complex conjugate of the matrix $F$. 

Let us recall that given an operator function $\Phi\in 
L^\infty(\T;E_*\tto E)$, the Toeplitz operator $T_\Phi: H^2(E_*)\to 
H^2(E)$ with symbol $\Phi$ is defined by 
$$
T_\Phi f := P_+(\Phi f), 
$$
where $P_+$ is the Riesz Projection (orthogonal projection onto 
$H^2$).

Considering the adjoint operator $(T_{\overline F})^* = T_{\overline 
F^*} =T_{F^T}$ one can conclude from here that $F$ is left invertible 
in $H^\infty$ if and only if the Toeplitz operator $T_{F^T}:H^2(E)\to 
H^2(E_*)$ is \emph{right} invertible. Since $F^T$ is an analytic 
function
$$
T_{F^T} f = F^T f, \qquad \forall f\in H^2(E), 
$$
and $F$ is left invertible in $H^\infty$ if and only if for any 
$g\in H^2(E_*)$ the equation
\begin{equation}
\label{0.1}
F^T f = g
\end{equation}
has a solution $g\in H^2(E)$ satisfying the uniform estimate
$
\|f\|_2 \le C\|g\|_2 
$.

The result that condition \eqref{C} implies (if $\dim 
E_*<+\infty$) left invertibility of the Toeplitz operator 
$T_{\overline F}$, or equivalently the solvability of the equation 
\eqref{0.1}, is called the \emph{Toeplitz Corona Theorem}. In the 
case of the unit disk $\D$ one can easily deduce the Matrix Corona 
Theorem  from the Toeplitz Corona Theorem by using the Commutant 
Lifting Theorem. 

The main result of this paper is  the Toeplitz Corona Theorem for the 
polydisk see Theorem \ref{t0.2} below.
To simplify the notation we used $F$ instead of $F^T$, so the 
condition \eqref{C} is replaced by the condition $FF^*\ge \delta^2 
I$. While in the polydisk it is not known how to get the Corona 
Theorem from the Toeplitz Corona Theorem (the Commutant Lifting Theorem 
for the polydisk is currently not known) the result seems to be of independent 
interest. In a particular case when $F$ from Theorem \ref{t0.2} is a 
row vector (a $1\times n$ matrix) this theorem was proved by K.~C.~Lin, see 
\cite{Lin1} or \cite{Lin2}. His approach involved using the Koszul complex to 
write down the $\pbar$-equations. Unfortunately, in several 
variables, unlike the one-dimensional case, higher order equations appear 
in addition to the $\pbar$-equation so the computation become quite 
messy. Moreover, it is not clear how to use his technique to get the 
result in the matrix case we are treating here since the Fuhrmann--Vasyunin 
trick of getting the matrix result from the result for a column (row) vector 
does not work to solve the Toeplitz Corona Theorem.

To prove the main result we use tools from complex differential 
geometry to solve $\db$-equations on  holomorphic vector bundles. In 
doing this we are following the ideas of M.~Andersson, see 
\cite{Andersson1} or \cite{Andersson2}, which in turn go back to B.~Berndtsson.

While our approach is quite similar to the one used by M.~Andersson, 
there are some essential differences.   To solve the $\db$-equation he uses a 
H\"ormander type approach with weights and a modification of a 
Bochner-Kodaira-Nakano-H\"ormander identity from complex geometry.  While our approach is more along the lines of T.~Wolff's proof and 
does not require anything more advanced than Green's formula.

We first use our technique to get an estimate in the Toeplitz Corona 
Theorem in the disk:

\begin{thm}
\label{t0.1}
Let $F \in H^\infty(\D;E\tto E_*)$, $\dim E_* =r<+\infty$, such that 
$\delta^{2}I\leq FF^{*}\leq I$ for some $0<\delta^{2}\leq\f{1}{e}$.  
For $1\leq p\leq\infty $ if $g\in H^p(\D;E_*)$ then  the equation
\begin{displaymath}
Ff=g
\end{displaymath}
has an analytic solution $f\in H^p(\D;E)$ with the estimate
\begin{equation}
\label{0.2}
\norm{f}_{p}\leq  \left( \f{C}{\delta^{r+1}}\log\f{1}{\delta^{2r}} 
+\frac1\delta\right)  \| g\|_p ,
\end{equation}
with $C= \sqrt{1+e^2} + \sqrt{e} + \sqrt{2} e \approx 8.38934$.
\end{thm}

For the $p=2$ case the above result with a different constant $C$ was obtained 
recently using a different method by T.~Trent \cite{Trent}.  The constant  he obtained was 
$C=2\sqrt e + 2\sqrt 2 e\approx 10.9859$. 

The result for all $p$ can be obtained from the case $p=2$ via the 
Commutant Lifting Theorem, but we present here a simple direct proof. 

\begin{rem*}Note, that we do not assume $\dim E<+\infty$ here. 
\end{rem*}

\

Using a simple modification of our proof in one dimension we are also 
able to get the following result in the polydisk:

\begin{thm}
\label{t0.2}
Let $F \in H^\infty(\D^n;E\to E_*)$, $\dim E_* =r<+\infty$, such that 
$\delta^{2}I\leq FF^{*}\leq I$ for 
some $0<\delta^{2}\leq\f{1}{e}$.  For $1<p<\infty $ if $g\in 
H^p(\D^n;E_*)$ then the equation
\begin{displaymath}
Ff=g
\end{displaymath}
has an analytic solution $f\in H^p(\D^n;E)$ with the estimate
\begin{equation}
\label{0.3}
\norm{f}_{p}\leq  
\left(\f{nCC(p)^n}{\delta^{r+1}}\log\f{1}{\delta^{2r}} + 
\frac1\delta\right)\| g\|_p ,
\end{equation}
where $C= \sqrt{1+e^2} + \sqrt{e} + \sqrt{2} e \approx 8.38934$, and 
$C(p)=1/\sin(\pi/p)$ the norm of the (scalar) Riesz projection from $L^p(\T)$ onto 
$H^p(\D)$.  For $p=2$ the estimate can be improved to
\begin{equation}
\label{0.4}
\norm{f}_2\leq\left( 
\f{\sqrt{n}C}{\delta^{r+1}}\log\f{1}{\delta^{2r}} +\frac1\delta 
\right) \norm{g}_2,
\end{equation}
with $C= \sqrt{1+e^2} + \sqrt{e} + \sqrt{2} e \approx 8.38934$.
\end{thm}

\subsection{Plan of the paper}
We will start with proving Theorem \ref{t0.1}  for $p=2$. 

In Section
\ref{s1} we set up the main estimate needed to prove the 
theorem. Section \ref{s2} is devoted to
a version of the Carleson Embedding Theorem and its analogue for 
functions defined on holomorphic vector
bundles, which will be later used to prove the main estimates. 

In Section \ref{s3} we peform computation of some derivatives and 
Laplacians that will be used in the
estimates. We also construct there subharmonic functions to be used 
in the embedding theorems. 

Section \ref{s4} deals with the main estimate for $p=2$; Section 
\ref{s5} explains how to use the
construction for other $p$. In Section \ref{s6} we treat the case of the
polydisk for $p=2$ and in Section
\ref{s7} we treat the case of general $p$. 

\section{Reduction to the main estimate}
\label{s1}

To prove Theorem \ref{t0.1} for $p=2$, for a given $g\in H^2:= 
H^2(E_*)$ with $\|g\|_2 = 1$, we need to
solve the equation
\begin{equation}
\label{1.1}
Ff = g, \qquad f \in H^2(E)
\end{equation}
with the estimate $\|f\|_2\le C=C(\delta, r)$.  By a normal families 
argument it is enough to suppose that
$F$ and
$g$ are  analytic in a neighborhood of $\D$.  Any estimate obtained 
in this case can be used to find an
estimate when $F$ is only analytic on $\D$.  Since $\delta^{2} I\leq 
FF^{*}\leq I$, it is easy to find a
non-analytic solution $f_0$ of \eqref{1.1}, 
$$
f_0 := \Phi g := F^*(FF^*)^{-1} g .  
$$
  
To make $f_0$ into an analytic solution, we need to find $v\in 
L^2(E)$ such that $f :=f_0 -v \in H^2$ and 
$v(z)\in \ker F(z)$ a.e.~on $\T$.   Then   
$$
Ff=F(f_{o}-v)=Ff_{o}-Fv=g,
$$
and we are done.  The standard way to find such $v$ is to solve a 
$\overline\partial$-equation with the 
condition $v(z)\in \ker F(z)$ insured by a clever algebraic trick. 
This trick also admits a ``scientific''
explanation, for one can get the desired formulas by writing a Koszul 
complex.  What we do in this paper
essentially amounts to solving the $\db$-equation $\overline\partial 
v = \overline\partial f_0$ on the
holomorphic vector bundle $\ker F(z)$. We mostly follow the  ideas of 
Matts Andersson found in
\cite{Andersson1}.  He used ideas from complex differential geometry 
to solve the corona problem by finding
solutions to the $\db$-equation on holomorphic vector bundles. 

Since our target audience consists of analysts, all  differential 
geometry will be well hidden. Our main 
technical tool will be Green's formula
\begin{equation}
\label{GF}
\int_{\T} u \, dm - u(0) = \frac{1}{2\pi} \int_{\D} \Delta u 
\log\frac{1}{|z|} dxdy.
\end{equation}
Instead of the usual Laplacian $\Delta= \frac{\p^2}{\p x^2} + 
\frac{\p^2}{\p y^2}$ it is more convenient for us to use the 
normalized one $\td :=\frac14\Delta = \db \p = \p \db$. If we denote 
by $\mu$ the measure defined by
$$
d\mu = \frac2\pi \log\frac{1}{|z|} dxdy, 
$$ 
then Green's formula can be rewritten as 
\begin{equation}
\label{GF1}
\int_{\T} u \, dm - u(0) =  \int_{\D} \td u \,d\mu. 
\end{equation}

\subsection{Set-up} To find the function $v$ we will use duality. We 
want $f_0- v \in H^2(E)$, therefore the equality
$$
\int_\T \La f_0, h \Ra \,dm = \int_\T \La v, h \Ra \,dm 
$$
must hold for all $h\in (H^2)^\perp$. Using Green's formula we get 
$$
\int_\T \La f_0, h \Ra \,dm = \int_\T \La \Phi g, h \Ra \,dm = 
\int_\D \p \db \left[ \La \Phi g, h \Ra \right] \,d\mu =
\int_\D \p  \left[ \La \db\Phi g, h \Ra \right] \,d\mu.
$$
Here we used the harmonic extension of $h$, so $h$ is anti-analytic 
and $h(0)=0$.  The functions $\Phi := F^*(FF^*)^{-1}$ and $g$ are 
already defined in the unit disk $\D$.

Now the critical moment: let $\Pi(z):=P\ci{\ker F(z)}$ be the 
orthogonal projection onto $\ker F(z)$, 
$\Pi=I-F^{*}\left(FF^{*}\right)^{-1}F$. Direct computation shows 
that  $\db \Phi= \Pi (\p \Phi)^* (FF^*)^{-1}$, so $\Pi \pbar \Phi = 
\pbar\Phi$. Therefore, if we define a vector valued function  $\xi$ 
on $\D$ by $\xi(z) := \Pi(z) h(z)$, then 
\begin{equation}
\label{1.3a}
\int_\D \p  \left[ \La \db\Phi g, h \Ra \right] \,d\mu = \int_\D \p  
\left[ \La \db\Phi g,\Pi h \Ra \right] \,d\mu = \int_\D \p  \left[ 
\La \db\Phi g, \xi \Ra \right] \,d\mu =: L(\xi). 
\end{equation}
Note, that $L$ is a conjugate linear functional, i.e.~$\overline L$ (defined by $\overline L(\xi):= \overline{L(\xi)})$ is a linear functional. 
Suppose we are able to prove the estimate 
\begin{equation}
\label{1.4}
| L(\xi) | \le C(r, \delta) \|\xi\|_2, \qquad \forall \xi = \Pi h, 
\quad h\in H^2(E)^\perp. 
\end{equation}
Then (by a Hilbert space version of the Hahn--Banach Theorem, which is 
trivial) $L$ can be extended to a bounded linear functional on 
$L^2(E)$, so there exists a function $v\in L^2(E)$, $\|v\|_2 \leq C$, 
such that 
$$
L(\xi) = \int_\T \La v, \xi \Ra \, dm , \qquad \forall \xi=\Pi h, 
\quad h\in H^2(E)^\perp. 
$$
Replacing $v$ by $\Pi v$ we can always assume without loss of 
generality that $v(z)\in \ker F(z)$ a.e.~on $\T$, so $Fv=0$.   By the construction
$$
\int_\T \La v, h \Ra \, dm = \int_\T \La v, \Pi h \Ra \, dm = L(\Pi 
h) = \int_\T \La \Phi g, h \Ra \, dm \qquad \forall h\in 
H^2(E)^\perp, 
$$
so $f:=f_0-v := \Phi g -v \in H^2(E)$ is the analytic solution we want 
to find.  It satisfies the estimate
$$
\|f\|_2 \le \|f_0\|_2 + \|v\|_2\le \frac1\delta \|g\|_2 + C(r, 
\delta) \|g\|_2. 
$$

Therefore, 
Theorem \ref{t0.1} would follow from the following proposition
\begin{prop}
\label{pr1.1}
Under the assumptions of Theorem \ref{t0.1}  the linear
functional
$L$ defined by
\eqref{1.3a} satisfies the estimate
$$
| L(\xi) | \le C(r, \delta) \|\xi\|_2, \qquad \forall \xi = \Pi h, 
\quad h\in H^2(E)^\perp, 
$$
with

$$
C(r, \delta) = \f{C}{\delta^{r+1}}\log\f{1}{\delta^{2r}}, 
$$
where $C= \sqrt{1+e^2} + \sqrt{e} + \sqrt{2} e $.
\end{prop}

In what follows we will need the following simple technical lemma 
that is
proved by direct computation. 

\begin{lm}
\label{l1.1}
For $\Pi$ and $\Phi$ defined above we have 
\begin{eqnarray*}
\p\Pi & = & -F^{*}\left(FF^{*}\right)^{-1} F'\Pi,\\
\pbar\Phi & = & \Pi\left( F'\right)^{*}\left(FF^{*}\right)^{-1},\\
\textnormal{and}\ \p\pbar\Phi & = & \p\Pi ( F')^* (FF^*)^{-1} - (\db 
\Phi) 
F' \Phi  = \p \Pi \db \Phi + 
(\p \Pi)^* \Phi  F' \Phi.
\end{eqnarray*}
\end{lm}

\begin{cor}
\label{c1.3}
For the projection $\Pi$ defined above we have 
$$
\Pi \partial \Pi =0, \qquad (\partial \Pi) \Pi =\p\Pi, \qquad (\db \Pi)\Pi
=0,
\qquad \Pi\db \Pi =\pbar\Pi.  
$$
\end{cor}

The above identities are well-known in complex differential geometry, but
we can easily get them from Lemma \ref{l1.1}. Namely, since $\Pi $ is the
orthogonal projection onto $\ker F$ we have $F\Pi=0$. Taking the adjoint we
get $\Pi F^*=0$ which implies $\Pi\partial \Pi =0$. The second identity
is trivial, and the last two are obtained from the first two by taking adjoints.

\section{Embedding theorems and Carleson measures}
\label{s2}

As is well known, Carleson measures play a prominent role in the 
proof 
of the Corona theorem, both in 
Carleson's original proof and in T.~Wolff's proof and subsequent 
modifications. It is also known to the
specialists, that essentially all
\footnote{By ``essentially all'' we mean here that a Carleson measure 
should first be mollified, to make it smooth, and then it can be 
obtained from the Laplacian of a subharmonic function. }
Carleson measures can be obtained from the Laplacian of a bounded 
subharmonic function. We will need the following well-known theorem, 
see \cite{Nikolskii}, which   was probably first proved by Uchiyama. 

\begin{thm}[Carleson Embedding Theorem]
\label{t2.1}
Let $\vf$ be a non-negative, bounded, subharmonic function. Then for 
any $f\in H^2(E)$
$$
\int_\D \td \vf (z) \|f(z)\|^2 d\mu(z) \le e \|\vf\|_\infty \|f\|_2^2.
$$
Here $d\mu = \frac{2}{\pi}\log\frac{1}{|z|}dxdy$, and $\td= \frac14 
\Delta = \p\db$.  
\end{thm}
\begin{proof}
Because of homogeneity, we can assume without loss of generality that 
$\|\vf\|_\infty =1$. Direct computation shows that 
$$
\td \bigl( e^{\vf (z)} \|f(z)\|^2\bigr) = e^\vf \td \vf \|f\|^2 + 
e^\vf \|\p \vf f + \p f \|^2 \ge  \td \vf \|f\|^2. 
$$
Then Green's formula implies 
$$
\int_\D \td \vf \, \|f\|^2\, d\mu \le \int_\D \td \bigl( e^\vf  
\|f\|^2\bigr)\, d\mu = 
\int_\T e^\vf  \|f\|^2 \, dm - e^{\vf(0)} \|f(0)\|^2 \le e\int_\T   
\|f\|^2 \, dm =e \, \|f\|_2^2.
$$
\end{proof}

\begin{rem}
It is easy to see, that the above Lemma implies the embedding 
$\int_\D \|f\|^2 \,d\mu \le C\int_\T \|f\|^2\, dm$ (with $C=e$) for 
all analytic functions $f$. Using the function $4/(2-\vf)$ instead of 
$e^\vf$ it is possible to get the embedding for harmonic functions 
with the constant $C=4$.  We suspect the constants $e$ and $4$ are 
the best possible for the analytic and harmonic embedding 
respectively. We cannot prove that, but it is known that $4$ is the 
best constant in the dyadic (martingale) Carleson Embedding Theorem. 
\end{rem}

We will need a similar embedding theorem for functions of form $\xi 
=\Pi h$, $h\in H^2(E)^\perp$. Such functions are not analytic or 
harmonic\footnote{To be precise, such functions are anti-holomorphic 
functions (with respect to the metric connection) on the holomorphic 
hermitian vector bundle $\ker F(z)$}, so the classical Carleson 
Embedding Theorem does not apply. As a result, the proof is more 
complicated, and the constant is significantly worse.

We will need several formulas.  Recall that  $\Pi(z) =P\ci{\ker 
F(z)}$ is the orthogonal projection onto $\ker F(z)$,  $\Pi = I - 
F^*(FF^*)^{-1}F$, and that $d\mu = \frac2\pi \log\frac{1}{|z|} dxdy$. 
\begin{lm}
\label{l2.3}
Let $\vf$ be a non-negative, bounded, subharmonic function in $\D$ 
satisfying
$$
\td \vf(z) \geq \|\p \Pi(z)\|^2, \qquad \forall z\in \D, 
$$
and let $K=\|\vf\|_\infty$.  Then for all $\xi$ of the form $\xi = \Pi 
h$, $h\in H^2(E)^\perp$
$$
\int_\D  \td \vf(z)\, \|\xi(z)\|^2 \, d\mu (z) \le e K e^K 
\|\xi\|_2^2 \, 
$$
and
$$
\int_\D  \|\db \xi \|^2\, d\mu \le (1 + eKe^K)\|\xi\|_2^2.
$$
\end{lm}
\begin{proof}
Let us take an arbitrary non-negative bounded subharmonic function 
$\vf$ and  compute $\td \left( e^\vf \|\xi\|^2\right)$. Corollary 
\ref{c1.3} implies that $\Pi \p \Pi =0$  and $\p \Pi \Pi = \p \Pi$. 
Therefore, using $\p h= 0$ we get  $\p \xi = \p \left( \Pi h \right) 
= \p \Pi h + \Pi \p h = \p \Pi h = \p \Pi \xi$, and so
$$
\La \p\xi , \xi\Ra = \La \p\xi , \Pi \xi\Ra = \La \p\Pi \xi , \Pi 
\xi\Ra = 0.
$$ 

Therefore 
$$
\p \left( e^\vf \|\xi\|^2\right) = e^\vf \p \vf \|\xi\|^2 + e^\vf \La 
\p \xi, \xi\Ra + e^\vf \La \xi, \db \xi\Ra = e^\vf \p \vf \|\xi\|^2 + 
e^\vf \La \xi, \db \xi\Ra. 
$$ 
Taking $\db$ of this equality (and again using $\La \xi ,\p\xi\Ra = 
0$) we get 
$$
\td \left( e^\vf \|\xi\|^2\right) = e^\vf \left( \td \vf \|\xi\|^2 + 
\|\db \vf \xi + \db \xi \|^2 + \La \xi, \td \xi \Ra \right). 
$$
To handle $\La \xi, \td \xi \Ra$ we take the $\p$ derivative of the 
equation $\La \xi , \p \xi \Ra =0$ to get 
$$
\La \p\xi, \p\xi \Ra + \La \xi, \db\p \xi \Ra =0, 
$$
and therefore $ \La \xi, \td \xi \Ra = - \| \p\xi \|^2 =- \|(\p \Pi) 
\xi \|^2$. 
Since $\vf\ge0$
\begin{align}
\notag
& \int_{\D}\left({\td \varphi}\|\xi\|^{2}-\|(\p\Pi) \xi\|^{2}\right) 
d\mu \le \\
\label{2.1}   &
 \qquad\qquad \int_{\D}\left({\td \varphi}\|\xi\|^{2} -\|(\p\Pi) 
\xi\|^{2} + \|\pbar\varphi\xi+\pbar\xi\|^{2} \right)e^{\varphi}d\mu = 
\int_{\T}e^{\varphi}\|\xi\|^{2}dm;
\end{align}
the  equality is just Green's formula (recall that $\xi(0)=0$).  In 
the last inequality replacing $\vf $ by $t\vf$, $t >1$ we get 
$$
\int_{\D}\left( t{\td\varphi}\|\xi\|^{2}-\|(\p\Pi) 
\xi\|^{2}\right)\,d\mu
\le 
\int_{\T}e^{t\varphi}\|\xi\|^{2}dm \le e^{tK} \|\xi \|_2^2.
$$
Now we use the inequality $\td \vf \geq \|\p \Pi\|^2$.  It implies 
$\td \vf \,\|\xi\|^2 - \|\p \Pi \xi \|^2 \geq 0$, and therefore
$$
(t-1) \int_{\D} {\td\varphi}\, \|\xi\|^{2}d\mu 
\le e^{tK} \|\xi\|_2^2 . 
$$
Hence
$$
\int_{\D} {\td\varphi}\, \|\xi\|^{2} d\mu \le
\min_{t> 1} \frac{e^{tK}}{t-1} \,  \|\xi\|_2^2 = e K e^{K} \|\xi\|_2^2
$$
(minimum is attained at $t=1+1/K$), and thus the first statement of 
the lemma is proved. 

To prove the second statement, put $\vf \equiv 0$ in \eqref{2.1} (we 
do not use any properties of $\vf$ except that $\vf \ge0$ in 
\eqref{2.1}) to get 
$$
\int_\D \left( \|\db \xi \|^2 - \|(\p \Pi) \xi \|^2 \right) \, d\mu  
= 
\int_\T \|\xi\|^2 \, dm = \|\xi\|_2^2. 
$$
But the second term can be estimated as
$$
\int_\D \|(\p \Pi) \xi \|^2 \, d\mu \le \int_\D \td \vf \, \| \xi 
\|^2 \, d\mu \le eK e^K\|\xi\|_{2}^{2},
$$
and therefore $\int_\D \|\db \xi \|^2 \, d\mu \le (1+ eK e^K) 
\|\xi\|_2^2$ .
\end{proof}

\section{Finding the correct subharmonic functions}
\label{s3}
There will be points in the proof where we would like to invoke 
Carleson's Embedding Theorem. To do so we will need a non-negative, 
bounded, subharmonic function.  In this section we construct the 
necessary subharmonic functions so they will be available when we 
finally estimate the integral in question.  With this in mind we 
define the two functions used and collect their relevant properties.  
First, we recall a basic fact that will aid in showing that the 
functions we construct are subharmonic.
\begin{lm}
\label{l3.1}
Let $A(t)$ be a  differentiable $n\times n$ matrix-valued function.  
Define the function $f(t)=\det(A(t))$.   Then
$$
f'(t)=\det(A(t))\tr(A^{-1}(t)A'(t)).
$$
\end{lm}
\begin{proof}
Fix a point $t$ and for brevity of notation let us use $A$ instead of
$A(t)$. Since $A(\cdot)$ is differentiable 
\begin{eqnarray*}
\det(A(t+h)) = \det(A+A'h+o(h)) & = & \det A\det(I+A^{-1}A'h+o(h))\\
 & = & \det A\prod(1+h\mu_{k}+o(h))
\end{eqnarray*}
where $\mu_{k}$ are the eigenvalues of $A^{-1}(t)A'(t)$.  Expanding this product
we 
have
$$
\prod(1+h\mu_{k}+o(h))=1+h\sum\mu_{k}+o(h)= 1+ h\tr (A^{-1}A') +o(h).
$$
Then
$$
\det(A(t+h)) =\det(A) +h\det(A)\tr(A^{-1}A')+o(h),
$$ 
which implies the desired formula for the derivative.
\end{proof}

Define the function $\varphi= \tr(\log(\delta^{-2}FF^{*})) 
=\log\left( \delta^{-2n}\det(FF^*) \right)$.  Then a straight forward 
application of the above lemma gives
\begin{eqnarray*}
\td\varphi & = & \p\db\varphi\\
 & = & \p[\tr((FF^{*})^{-1}F(F')^{*})]\\
 & = & \tr[(FF^{*})^{-1}F'\Pi(F')^{*}]
\end{eqnarray*}
with the last line following by substitution of $\Pi$.  For another 
approach to this computation see \cite{Trent}.  Using the identities 
$\Pi^2 = \Pi$, $\tr(AB)=\tr (BA)$, and recalling that 
$$
\p \Pi = -F (FF^*)^{-1} F' \Pi 
$$ 
we get
\begin{eqnarray*}
\td \varphi & = & 
\tr\left[\left(FF^{*}\right)^{-1}F'\Pi(F')^{*}\right]\\
 & = & 
\tr\left[F^{*}\left(FF^{*}\right)^{-1}F'\Pi\Pi(F')^{*}\left(FF^{*}\right)^{-1}F\right]\\ 
& = & \tr\left[\p\Pi(\p\Pi)^{*}\right]\\
 & \geq & \norm{\p\Pi}^{2},
\end{eqnarray*}
with the last inequality following since 
$\tr[AA^{*}]\geq\norm{A}^{2}$.  This function will play a prominent 
role in the estimation of certain integrals.  We should also note 
that 
$$
0\leq\varphi\leq K:= \log\f{1}{\delta^{2n}}.
$$  

We will also need another function to help in the estimation of the 
linear functional $L$ in question.  Let $\lambda=\tr(\left(FF^{*}\right)^{-1})$.  
A simple  computation 
gives, 
\begin{eqnarray*}
\td \lambda & = & \tr\left[\Phi^* (F')^{*}\left(FF^{*}\right)^{-1}F' 
\Phi \right]\\
 &  & 
-\tr\left[\left(FF^{*}\right)^{-1}F'\Pi(F')^{*}\left(FF^{*}\right)^{-1}\right]\\ 
& \geq & \tr\left[\Phi^* (F')^{*}\left(FF^{*}\right)^{-1}F'\Phi^* 
\right]-\delta^{-2}\tr\left[\p\Pi(\p\Pi)^{*}\right].
\end{eqnarray*}
Now we define the function $\psi=\lambda+\delta^{-2}\varphi$. Then, 
recalling that $\Phi = F^*(FF^*)^{-1}$ we get
\begin{eqnarray*}
\td \psi & \geq & \tr\left[\Phi^*(F')^{*}\left(FF^{*}\right)^{-1}F' 
\Phi \right] 
\\
 & = & \tr\left[\Phi F'\Phi(\Phi F'\Phi)^{*}\right]\\
 & \geq & \left\| \Phi F\Phi \right\|^{2} .
\end{eqnarray*}
So $\psi$ is subharmonic and 
$0\leq\psi\leq\f{n}{\delta^{2}}+\f{1}{\delta^{2}}\log\f{1}{\delta^{2n}}$.  
We should note that the assumption $0<\delta^{2}\leq\f{1}{e}$ implies 
$\log\delta^{-2}\geq 1$. This gives 
$$
0\leq\psi  \le L :=  \f{2}{\delta^{2}}\log\f{1}{\delta^{2n}}.
$$

\section{Estimating the integral}
\label{s4}
Now we need to estimate $L(\xi)$.  Computing $\p$ of the inner 
product  we get
\begin{eqnarray*}
L(\xi) & = & \int_{\D}\p \left[ \La \pbar \Phi g,\xi\Ra \right] d\mu\\
 & = & \int_{\D}\La\p\pbar\Phi g,\xi\Ra d\mu+  \int_{\D}\La\pbar\Phi 
g',\xi\Ra d\mu+\int_{\D}\La\pbar\Phi g,\pbar\xi\Ra d\mu
\\ & = & I+II+III .
\end{eqnarray*}
We need to estimate each of the above integrals as closely as 
possible.  Each integral has a term involving derivatives of $\Pi$, 
$g$ and $\xi$. The idea is to separate the integrals using 
Cauchy-Schwarz, giving one derivative to each term.  

We now estimate the first integral.  Recalling that $\p\pbar \Phi = 
\p \Pi \pbar \Phi +(\p \Pi)^* \Phi F' \Phi$ we get 
$$
I = \int_{\D}\La\p\pbar\Phi g,\xi\Ra d\mu   
=  \int_{\D}\left\{ \La\p\Pi\pbar \Phi g,\xi\Ra + \La(\p\Pi)^{*}\Phi 
F' \Phi g,\xi\Ra \right\} d\mu .
$$
Since $(\p \Pi)^* \Pi =0$ we have $(\p \Pi)^* \xi =0$, and so 
$\La\p\Pi\pbar \Phi g,\xi\Ra =0$. Therefore
$$
I= \int_{\D} \La(\p\Pi)^{*}\Phi F' \Phi g,\xi\Ra  d\mu  =
\int_{\D} \La \Phi F' \Phi g,  (\p \Pi)\xi\Ra  d\mu , 
$$
and the Cauchy-Schwarz inequality implies
$$
|I| \le \left( \int_\D \| \Phi F' \Phi g \|^2 d\mu \right)^{1/2} 
\left( \int_\D \| (\p \Pi)\xi \|^2 d\mu \right)^{1/2}.
$$

To estimate the second factor we use Lemma \ref{l2.3}. Recall that 
the function 
$$
\vf= \log\left( \delta^{-2n} \det(  FF^*) \right), 
$$ 
constructed in Section \ref{s3} satisfies the inequalities
\begin{equation}
\label{4.1}
\td \vf \ge \|\p \Pi\|^2, \qquad\text{and} \qquad 0\le \vf \le K := 
\log\delta^{-2n}. 
\end{equation}
Therefore, Lemma \ref{l2.3} implies 
$$
\int_\D \| (\p \Pi)\xi \|^2 d\mu \le e K e^K \|\xi\|_2^2 = e 
\delta^{-2n} \log\delta^{-2n} \|\xi\|_2^2 .
$$

To estimate the first factor, notice that the function $\psi$ 
constructed in Section \ref{s3} satisfies
$$
\td \psi \ge \| \Phi F' \Phi \|^2, \qquad \text{and}\qquad 0\le \psi 
\le L:= 2\delta^{-2}\log \delta^{-2n}. 
$$
Then the Carleson Embedding Theorem (Theorem \ref{t2.1}) implies
$$
\int_\D \| \Phi F' \Phi g \|^2 d\mu \le e L \|g\|_2^2 =  
2e\delta^{-2}\log \delta^{-2n} \|g\|_2^2, 
$$
and thus
$$
|I| \le \sqrt{KL}\,\|\xi\|_2 \|g\|_2 = \frac{\sqrt2 e}{\delta^{n+1}} 
\log\delta^{-2n}  \|\xi\|_2  \|g\|_2. 
$$

Now we estimate $II$.  By the Cauchy-Schwarz inequality, we have 
\begin{eqnarray*}
| II |  & \leq & \int_{\D} | \La\pbar\Phi g',\xi\Ra |  d\mu\\
 & \leq & \left(\int_{\D}\norm{\pbar\Phi}^{2}\|\xi\|^{2} 
d\mu\right)^{1/2}\left(\int_{\D}\|g'\|^{2} d\mu\right)^{1/2}.
\end{eqnarray*}
Observe that $\td \|g\|^{2}= \|g'\|^2$ since $g$ is holomorphic. So, 
applying  Green's Theorem to the second factor we get
$$
\int_{\D}\|g'\|^{2} d\mu = \int_\T \|g\|^2 d m - \|g(0)\|^2 \le 
\|g\|_2^2. 
$$
To estimate the first integral, notice, that  
$$
\|\Phi \|^2 = \| \Phi^* \Phi \| = \| (FF^*)^{-1}\| \le \delta^{-2}
$$
(recall that $\Phi = F^*(FF^*)^{-1}$).  Since $\pbar \Phi = -(\p 
\Pi)^*\Phi$, we can estimate
$$
\| \pbar \Phi \|^2 = \| ( \pbar \Phi )^* \pbar \Phi \| = \| \Phi^* \p 
\Pi (\p \Pi)^* \Phi \| 
\le \| \p \Pi (\p \Pi)^*\| \cdot \| \Phi \|^2 
 \le \delta^{-2}  \|\p\Pi \|^2. 
$$
Therefore (see \eqref{4.1}), $\|\pbar \Phi\|^2 \le \delta^{-2} \td 
\vf$, where $\vf  =  \log\left( \delta^{-2n} \det(  FF^*) \right)$ is 
the subharmonic function constructed in Section
\ref{s3}. Applying Lemma \ref{l2.3} we get
$$
\int_{\D}\norm{\pbar\Phi}^{2}\|\xi\|^{2} d\mu \le \delta^{-2} 
\int_{\D}\td \vf \|\xi\|^{2} d\mu \le \delta^{-2} e K e^K 
\|\xi\|_2^2, 
$$
where $K= \log\delta^{-2n}$, see \eqref{4.1}. Joining the estimates 
together, we get 
$$
|II | \le \delta^{-1} \sqrt{eK} e^{K/2} \|g\|_2 \|\xi\|_2 \le 
\delta^{-1}
\sqrt{e} K e^{K/2} \|g\|_2 \|\xi\|_2 = 
\frac{\sqrt{e}}{\delta^{n+1}}\log\delta^{-2n}   \|g\|_2 \|\xi\|_2
$$
(since $\delta^{2}\leq\f{1}{e}$, the value of $K$ satisfies 
$K^{1/2}\leq K$).

Finally moving on to integral $III$.  Using Cauchy-Schwarz, we have
\begin{eqnarray*}
| III | & \leq & \int_{\D}\abs{\La\pbar\Phi g,\pbar\xi\Ra} d\mu\\
 & \leq & \left(\int_{\D}\norm{\pbar\Phi}^{2}\|g\|^{2} 
d\mu\right)^{1/2}\left(\int_{\D}\|\pbar\xi\|^{2} d\mu\right)^{1/2}.
\end{eqnarray*}
As we already have shown above, $\|\pbar \Phi \|^2 \le \delta^{-2} 
\td \vf$.  The Carleson Embedding Theorem (Theorem \ref{t2.1}) implies
$$
\int_\D \| \pbar \Phi  \|^2 \| g\|^2 d\mu \le \delta^{-2} \int_\D \td 
\vf \|g\|^2 d\mu 
\le \delta^{-2} e K \| g \|_2^2 .
$$
Using Lemma \ref{l2.3} we can estimate
$$
\int_\D \| \pbar \xi \|^2 d\mu \le (1 + e K e^K) \|\xi\|_2^2 \le 
(e^{-1} + e) K e^K \|\xi \|_2^2.
$$
Here we are using the fact that $K\ge1$ for $\delta^2\le 1/e$. 
Combining the estimates, we get 
$$
| III | \le \sqrt{1+e^2} K e^{K/2} \|g\|_2 \|\xi \|_2 \ =    \frac{ 
\sqrt{1+e^2}}{\delta^{n+1}} \log\delta^{-2n}\|g\|_2 \|\xi\|_2. 
$$

Joining the estimates for $I$, $II$, $III$ we get 
\begin{prop}
\label{pr4.1}
Under the assumptions of Theorem \ref{t0.1}  the linear
functional
$L$ defined by
\eqref{1.3a} satisfies the estimate
$$
| L(\xi) | \le C(r, \delta) \|\xi\|_2, \qquad \forall \xi = \Pi h, 
\quad h\in H^2(E)^\perp, 
$$
with 
$$
C(r, \delta) = \f{C}{\delta^{r+1}}\log\f{1}{\delta^{2r}}, 
$$
where $C= \sqrt{1+e^2} + \sqrt{e} + \sqrt{2} e $.

\end{prop}
Proposition \ref{pr4.1} is just a restatement of Proposition \ref{pr1.1}, and this then proves Theorem \ref{t0.1} for the case of $p=2$.  Note, that the constant $C$ is a bit better than the constant $2\sqrt2 e + 
2\sqrt{e} \approx 10.9859$ obtained by T.~Trent in \cite{Trent}.

\section{The $H^{p}$ Corona Problem in the Disk}
\label{s5}

Now we indicate how we can use the $H^{2}$ result to figure out the 
$H^{p}$ result.  We can use much of 
the same approach as in the $H^{2}(E)$ case.  Our goal is to solve 
the equation
$$
Ff=g,\qquad f\in H^{p}(E)
$$
for the given $g\in H^{p}(E_*)$, with $\norm{g}_p=1$, and furthermore 
we want the estimate $\norm{f}_p\leq C$.  Again we will have the 
obvious non-analytic solution to the problem 
$$
f_0:=\Phi g:=F^*(FF^{*})^{-1} g.
$$
To make this into an analytic solution we will need to find a 
function $v\in L^{p}(E)$ such that $f_0-v\in H^{p}$ and $v(z)\in\ker 
F(z)$.  This will be accomplished by duality.  As in the $H^{2}(E)$ 
case we need
$$
\int_{\T}\La f_0, h\Ra dm = \int_{\T}\La v,h\Ra dm
$$
to hold for all $h\in H^{p}(E)^{\perp}=H^{q}_0(E)$ (this uses the 
standard duality of $H^{p}$ spaces see \cite{Garnett} or 
\cite{Nikolskii}).  Again we can ensure that $v\in \ker F(z)$ since 
$\db\Phi=\Pi\db\Phi$.  So we need to get an estimate on the linear 
functional
$$
L(\xi)=\int_{\D}\p\left[\La\db\Phi g,\xi\Ra\right]d\mu
$$
with $\xi=\Pi h$ and $h\in H^{p}(E)^{\perp}$.  If we can then prove 
that
$$
\abs{L(\xi)}\leq C\norm{\xi}_{q}
$$
then by  the Hahn-Banach Theorem and 
duality in $L^{p}$ spaces with values in a Hilbert space we would have the
existence of a function 
$v\in L^{p}(E)$ with $\norm{v}_p\leq C$, such that 
$$
L(\xi)=\int_{\T}\La v,\xi\Ra dm,\quad\forall\xi=\Pi h,\quad h\in 
H^{p}(E)^{\perp}.
$$
Then replacing $v$ by $\Pi v$ we can assume without loss of 
generality that $v(z)\in \ker F(z)$ a.e. on $\T$.  But then the 
construction would give,
$$
\int_{\T}\La v,h\Ra dm=\int_{\T}\La v,\Pi h\Ra dm = L(\Pi 
h)=\int_{\T}\La\Phi g,h\Ra dm,\quad\forall h\in H^{2}(E)^{\perp},
$$
so $v-f_0\in H^{p}(E)$.  So we only need to show how to prove the 
estimate
$$
\abs{L(\xi)}\leq C\norm{\xi}_q.
$$

To prove this estimate we will have to consider the cases $1<p<2$, 
$2<p<\infty$, $p=1$, and $p=\infty$.  The idea is to multiply the 
function $g$ and $\xi$ by some scalar function that has the correct norm, 
but places each of these functions in $H^{2}(E)$.  

First look at the case $1<p<2$.  In this case we will need to 
multiply by a scalar valued outer function that has the same norm as 
$g(z)$, namely we will have
$$
g_{out}(z)=\abs{g(z)}\ \textnormal{a.e. on }\T.
$$
Since we need to multiply by an outer function it will be important 
that $\log\norm{g}\in L^{1}(\D)$.  But since $g$ is holomorphic this 
is no problem.  So we the multiply $g$ by $g_{out}^{p/2 -1}$ and 
multiply $\xi$ by $\overline{g}_{out}^{1-p/2}$ since we need to take 
into consideration the conjugate.  Then call

\begin{eqnarray*}
\tilde{g}(z) & = & g_{out}^{p/2-1}(z)g(z)\ \textnormal{and}\\
\tilde{\xi}(z) & = & \overline{g}_{out}^{1-p/2}(z)\xi(z).
\end{eqnarray*}

Then computation using H\"older's Inequality with $q$ the conjugate 
exponent to $p$ gives $\norm{\tilde{g}}_2=\norm{g}_p^{p/2}$ and 
$\norm{\tilde{\xi}}_2\leq\norm{\xi}_q\norm{g}_p^{\f{p^2}{2(2-p)}}$.

Now using the $H^{2}$ result just proved we have that 
\begin{eqnarray*}
\abs{L(\xi)}=\abs{L(\tilde{\xi})} & \leq & 
\f{C}{\delta^{n+1}}\log\delta^{-2n}\norm{\tilde{g}}_2\norm{\tilde{\xi}}_2\\ 
& \leq & \f{C}{\delta^{n+1}}\log\delta^{-2n}\norm{g}_p\norm{\xi}_q.
\end{eqnarray*}

But this is the result that we sought to prove, and so we have taken 
care of the case $1<p<2$.

To deal with the case when $2<p<\infty$ is analogous, except in this 
case we need to multiply by an outer function corresponding to 
$\xi$.  The manipulations are identical to above, but the existence 
of the outer function requires a brief explanation.  We need to make 
sure that we can find an outer function that has the same norm as 
$\xi=\Pi h$.  Now recall that $\Pi = I-F^*(FF^{*})^{-1}F$.  Since $F$ 
is analytic, we may assume that $\Pi$ is real analytic in $\D$.  By 
density we may assume that $h$ is in fact a polynomial, and as such 
has zeros of at most finite multiplicity in $\D$ and a.e. on $\T$.  
Then $\xi=\Pi h$ will have zeros of finite multiplicity and we will 
be able to construct $\xi_{out}$.

So it only remains to deal with the cases $p=1$ and $p=\infty$.  In 
the $p=\infty$ case the above argument works the same way.  We need 
to show that $\abs{L(\xi)}\leq C\norm{g}_{\infty}\norm{\xi}_1$.  Then 
we may use the fact that the dual of $L^1$ is $L^\infty$, and so the 
linear functional $L(\xi)$ can then be identified by a pairing with a 
$L^{\infty}(E)$ function.  Again, the idea will be to multiply the 
function $\xi$ by an appropriate scalar outer function with the 
correct norm.  Let $\xi_{out}(z)=\abs{\xi(z)}$ a.e. on $\T$.  Then 
define
\begin{eqnarray*}
\widetilde{g}(z) & = & \overline{\xi}_{out}^{1/2}(z)g(z)\\
\widetilde{\xi}(z) & = & \xi_{out}^{-1/2}(z)\xi(z).
\end{eqnarray*}
Computation then yields that 
$\norm{\widetilde{\xi}}_{2}=\norm{\xi}_{1}^{1/2}$, and 
$\norm{\widetilde{g}}_{2}\leq\norm{g}_{\infty}\norm{\xi}_{1}^{1/2}$, 
and combining this with the $H^{2}$ result we have that 
$\abs{L(\xi)}\leq C\norm{g}_\infty\norm{\xi}_1$, proving the 
$H^{\infty}$ result.

The case $p=1$ requires just a little more work since $L^{1}$ is not 
the dual of $L^{\infty}$.  In this case we need to show 
$\abs{L(\xi)}\leq C\norm{g}_1\norm{\xi}_{\infty}$.  We begin by 
multiplying by a suitable outer function.  Define
\begin{eqnarray*}
\widetilde{g}(z) & = & g_{out}^{-1/2}(z)g(z)\\
\widetilde{\xi}(z) & = & \overline{g}_{out}^{-1/2}(z)\xi(z).
\end{eqnarray*}

Then we have $\norm{\widetilde{g}}_2 =\norm{g}_1^{1/2}$ and 
$\norm{\widetilde{\xi}}_2\leq\norm{g}_{1}^{1/2}\norm{\xi}_\infty$.  
So we have $\abs{L(\xi)}\leq C\norm{g}_1\norm{\xi}_{\infty}$, as 
needed.  Now this implies that the functional $L$ can be extended to 
a linear functional on $L^{\infty}(E)$.  Since the continuous 
functions on the unit circle are contained in $L^{\infty}(E)$ there 
will exist a vector-valued measure $\nu$ such that 
$$
L(\xi)=\int_{\T}\xi d\nu.
$$
Without loss of generality replace $\nu$ with $\Pi\nu$, then
$$
\int_{\T}h d\nu=\int_{\T}\xi d\nu = L(\Pi h)=\int_{\T}\La f_0, h\Ra 
dm.
$$
Then re-writing this, and treating $f_0 dm$ as a vector valued 
measure we have
$$
\int_{\T}h(f_0 dm-d\nu)=0\qquad\forall h\in 
H^{1}(E)^{\perp}=H_0^{\infty}(E).
$$
So in particular we have that the measure annihilates all analytic 
polynomials.  Then applying the F. $\&$ M. Riesz Brothers Theorem, see \cite{Nikolskii}, we can conclude that the measure $f_0 dm-d\nu$ 
is absolutely continuous with respect to Lebesgue measure, and 
moreover it is an analytic measure meaning $f_0 dm-d\nu=(f_0 -v)dm$ 
with $v\in H^{1}(E)$.  Then this implies that 
$$
\int_{\T}\La f_0-v,h\Ra dm=0\qquad\forall h\in H^{1}(E)^{\perp},
$$
or $f_0-v\in H^{1}(E)$, proving the $p=1$ case.

\section{The $H^{2}$ Corona Problem in the Polydisk}
\label{s6}
In the following sections we will be considering operator-valued functions that 
are taking values from the polydisk $\D^{n}$ to appropriate Hardy 
classes.  We begin with the $H^{2}(E)$ case.  The general goal 
from previous sections has not changed.  We want, for a given $g\in 
H^2:= H^2(\D^n;E_*)$ with $\|g\|_2 = 1$, to solve the equation
\begin{equation}
\label{6.1}
Ff = g, \qquad f \in H^2(\D^n;E)
\end{equation}
with the estimate $\|f\|_2\le C$.  Again by a normal families 
argument it is enough to suppose that $F$ and $g$ are analytic in a 
neighborhood of $\D^{n}$ because any estimate obtained can be used to 
get an estimate when $F$ is only analytic in $\D^n$.  It is still 
easy to find a non-analytic solution $f_0$ of \eqref{6.1}, 
$$
f_0 := \Phi g := F^*(FF^*)^{-1} g,   
$$
because we have $\delta^{2} I\leq FF^{*}\leq I$.  We will again need 
to find a $v\in L^{2}(\T^n;E)$ such that $f:=f_0-v\in H^{2}(\D^n;E)$ 
with $v(z)\in\ker F(z)$ a.e.~on $\T^{n}$.  Our approach is straight 
forward reduction to the one variable case, unfortunately this 
approach will not yield a proof of the $H^{\infty}$ Corona problem on the polydisk
since the projections are not bounded when $p=\infty$.

We will denote a point in $\D^n$ or $\T^n$ by $\bz=(z_1, z_2, \ldots, z_n)$. We
will use the symbol $\bz_j$ for $\bz$ without the coordinate $z_j$ and, slightly
abusing notation, we can then write 
$\bz=(\bz_j, z_j)=(z_j,\bz_j)$.

Let $H^p_j=H^p_j(\D^n;E)$ be a subspace of $L^p(\T^n;E)$ consisting of all functions
analytic in $z_j$, i.e.
\begin{equation}
\label{6.1.1}
H^p_j(\D^n;E) :=\left\{ f\in L^p(\T^n, E) : f(\bz_j, \cdot)\in H^p (\D;E) \text{
for almost all } \bz_j\in\T^{n-1} \right\}.
\end{equation}

\subsection{Lemmas about decompositions}

\begin{lm}
\label{l6.1}
Any $h\in H^{2}(\D^n;E)^\perp$ can be written as 
$h=\sum_{j=1}^{n}h_{j}$ with 
$h_j\in H^2_j(\D^n;E)^\perp$. 
\end{lm}

\begin{proof}
Let $P_j :=P\ci{H^2_j}$ be the orthogonal projection onto $H^2_j: = H^2_j(\D^n;E)$.
We can decompose $h$ in the following way,
$$
h= P_1 h + (I-P_1)h = h^1 + h_1 \qquad h_1\in H^2_1(\D^n;E), \ h^1 = P_1h.
$$
Similarly, 
$$
h^1= P_2 h^1 + (I-P_2)h^1 = h^2 + h_2 \qquad h_2\in H^2_2(\D^n;E), \ h^2 =
P_2P_1h.
$$
Continuing the procedure we get
$$
h^{k-1}= P_k h^{k-1} + (I-P_k)h^{k-1} = h^k + h_k \qquad h_k\in H^2_k(\D^n;E), \
h^k = P_k \cdots P_2P_1h.
$$
Combining everything we get 
$$
h=h_1 + h_2 +\ldots + h_n + h^n, \qquad h^n=P_nP_{n-1}\cdots P_1h
$$
which proves the lemma, because the assumption $h\in H^2(\D^n; E)^\perp$ implies
that $h^n= P_n\ldots P_2 P_1 h =0$
\end{proof}

We also are going to need an analogue of Lemma \ref{l6.1} dealing with the decomposition of 
functions on the holomorphic vector bundle $\Pi H^2$, i.e.~for the functions of the form $\xi=\Pi h$, $h\in 
H^{2}(\D^{n};E)^{\perp}$. To state this lemma we need some auxiliary 
definitions. Let 
\begin{equation}
\label{6.2a}
K(\D^{n};E):=\textnormal{clos}(\Pi (H^{2}(\D^{n})^{\perp})),
\end{equation}
and
\begin{equation}
\label{6.3a}
K_j(\D^{n};E):=\textnormal{clos}(\Pi (H_{j}^{2}(\D^{n})^{\perp})),\quad\forall\ 
j=1,\ldots,n,
\end{equation}
%%
%where $H^{2}_{j}(\D^{n})$ is as in Lemma \ref{l6.1}.  
%We will also let $P_{K_j}$ denote the orthogonal projection in $\Pi L^2=\Pi L^2(\D^n;E)$ 
%onto the subspace $K_j$ and $Q_{K_j}:=I-P_{K_j}$ be the orthogonal projection onto $K_j^{\perp}$.

\begin{lm}
\label{l6.2}
Let $\xi\in K$, then $\xi=\sum_{j=1}^{n}\xi_{j}$ with $\xi_j\in K_j$ 
for $j=1,
\ldots,n$ and 
$$\norm{\xi}^{2}_{2}=\sum_{j=1}^{n}\norm{\xi_j}_{2}^{2}.$$
\end{lm}

To prove Lemma \ref{l6.2} we will need a few other lemmas.  The first one is a simple fact about the geometry of a Hilbert space.

\begin{lm}
\label{l6.2.1}
Let $X$ be a subspace of a Hilbert space $H$, and let $\Pi$ be some orthogonal projection in $H$.
Then $\ran \Pi =\Pi H$ is decomposed into the orthogonal sum
$$
\Pi H = \clos( \Pi X )\oplus (X^\perp\cap \Pi H). 
$$ 
\end{lm}
\begin{proof}
The proof is a simple exercise in functional analysis,  and we leave it to the reader. 
\end{proof}

Define the subspaces
\begin{equation}
\label{6.5}
Q(\D^{n};E):=H^2\cap \Pi L^2, \qquad Q_j(\D^{n};E) := H^2_j \cap \Pi L^2. 
\end{equation}
Applying the above lemma to $H= L^2$ and $X=(H^2)^\perp$ or $X=( H^2_j)^\perp$ we get the following result.
\begin{cor}
\label{c6.4}
The subspace $\Pi L^2 =\Pi L^2(\D^n;E)$, $n=1, 2, 3, \ldots$ admits the orthogonal decompositions
$$
\Pi L^2 = K\oplus Q, \qquad\Pi L^2 = K_j\oplus Q_j,
$$
with the subspaces $K:= K(\D^n;E)$, $K_j:=K_j(\D^n;E) $, $Q:=Q(\D^{n};E)$ and $Q_j:=Q_j(\D^{n};E)$ defined by \eqref{6.2a}, \eqref{6.3a} and \eqref{6.5} respectively. 
\end{cor}

\begin{rem}
\label{r6.5}
Note, that the orthogonal projections $P\ci{K_j}$ and $P\ci{Q_j}$ are essentially ``one-variable'' operators. 
Namely, to perform the projection $P_{Q_j}$ on the function $\xi\in \Pi L^2$ we simply need  to perform for each $\bz_j\in \T^{n-1}$ (recall that $\bz =(z_j, \bz_j)$) the ``one-variable'' projection $P\ci{Q(\bz_j)}$ onto the subspace 
$$
Q(\bz_j) := H^2(\D;E)\cap \Pi(\,\cdot\,,\bz_j) L^2(\D;E)\subset H^2 =H^2(\D; E), 
$$
and similarly for the projection $P_{K_j}$.

Indeed, if 
$$
\xi^1 (\,\cdot\,, \bz_j) := P_{Q(\bz_j) }\xi(\,\cdot\,, \bz_j) \qquad \text{for almost all }\bz_j\in \T^{n-1},
$$ 
then clearly 
$$
\xi^1(\,\cdot\,, \bz_j) \in H^2 (\D;E)\cap \Pi(\,\cdot\,, \bz_j) L^2(\T)\qquad \text{for almost all }\bz_j\in \T^{n-1},
$$ 
so $\xi^1\in H^2(\D^n;E) \cap \Pi L^2(\D^n;E)$. 
Moreover, for $\xi_1:=\xi-\xi^1$ and any $\eta\in H^2(\D^n;E)\cap \Pi L^2$
$$
\int_\T \Bigl\La\xi_1(z_j, \bz_j), \eta(z_j,\bz_j) \Bigr\Ra dm(z_j)=0 \qquad \text{for almost all } \bz_j\in \T^{n-1}, 
$$
and integrating over other variables $\bz_k$ we get that $\xi_1\perp \eta$. 
\end{rem}

The following two lemmas says that in many respects the projection $P\ci{Q_j}$ behaves like the projection $I-P_j$ from Lemma \ref{l6.1}. 
\begin{lm}
\label{l6.3}
Let $H^2=H^2(\D^2;E)$ and let  $Q$ and $Q_{j}$, $j=1, 2$, be the subspaces as defined above in \eqref{6.5}. 
 Then for the orthogonal projections $P\ci{Q_j}$ onto the subspaces $Q_j$ we have 
$$
P_{Q_1}P_{Q_2}= P_{Q_2}P_{Q_1}=P_{Q} \,.
$$
\end{lm}
\begin{proof}  
It follows from the definition of $Q$ and $Q_j$ and from the inclusion $H^2\subset H^2_j$ that 
$$
Q= \Pi L^2 \cap H^2 \subset \Pi L^2 \cap H^2_j = Q_j 
$$
we can conclude that for $\xi\in Q$ we have $P\ci{Q_j} \xi = \xi$, $j=1, 2$. 

Since by Corollary \ref{c6.4} we have the orthogonal decomposition $\Pi L^2 = K\oplus Q$, 
 to prove the lemma we need to show that   the equalities $P\ci{Q_2}P\ci{Q_1}\xi =0$, $P\ci{Q_1}P\ci{Q_2}\xi =0$ hold for all $\xi \in K$. Clearly, it is sufficient to prove only one, say the first as the second can be obtained by interchanging indices. 

Consider the orthogonal decomposition of $\xi\in K$,
$$
\xi =P\ci{K_1} \xi + P\ci{Q_1} \xi =: \xi_1 + \xi^1. 
$$
To prove that $P_{Q_2}P_{Q_1}\xi=0$ we need to show that $\xi^{1}\in K_2$.

By definition $\xi^1\perp K_1 := \textnormal{clos}(\Pi((H^2_1)^\perp))$, and since $\Pi((H^2_1)^\perp) \supset (H^2_1)^\perp \cap \Pi L^2$, we can conclude that 
$$
\xi^1\perp (H^2_1)^\perp \cap \Pi L^2.
$$
We know that $\xi, \xi_1\in K$ ($\xi_1\in K$ because $K_1\subset K$), so $\xi^1\in K$. By Corollary \ref{c6.4}, 
$$
\xi^1\perp Q:= H^2 \cap \Pi L^2.
$$
Combining the above two orthogonality relations we get 
$$
\xi^1 \perp \left((H_1^2)^{\perp}+H^2\right)\cap\Pi L^{2},
$$
and since in the bidisk $H_2^2 \subset (H^2_1)^\perp + H^2$, we get that 
$$
\xi^1\perp  \Pi L^2 \cap H^2_2 =: Q_2
$$
i.e.~that $\xi^1 \in K_2$. 
\end{proof}

As an important corollary we get the following lemma. 

\begin{lm}
\label{l6.4}
On $\Pi L^2: =\Pi L^{2}(\T^{n};E)$ we have
$$
P_{Q_k}P_{Q_j}= P_{Q_j}P_{Q_k}= P_{Q_k\cap Q_j}= P_{H_{jk}^{2}\cap\Pi L^{2}}\quad\forall 1\leq k,j\leq n,
$$
where $H_{jk}^{2}(\D^{n}):=H_{j}^{2}(\D^{n})\cap H_{k}^{2}(\D^{n})$.  Furthermore, this implies
$$
P_{Q_1}\ldots P_{Q_n}= P_{Q_n}\ldots P_{Q_1}=P_{H^{2}\cap\Pi L^{2}}.
$$
\end{lm}

One can think of the space $H_{jk}^{2}(\D^{n})$ as the space of functions in $L^{2}(\T^{n})$ which are, upon fixing the other variables, holomorphic in both the $j$th and $k$th variable.

\begin{proof}
The first part of the lemma follows immediately from Lemma \ref{l6.3}, because we can just ``freeze'' all variables except $z_j$ and $z_k$. 
Namely, to perform the projection $P_{Q_j}$ on the function $\xi\in \Pi L^2$ we simply need  to perform for each $\bz_j\in \T^{n-1}$ (recall that $\bz =(z_j, \bz_j)$) the ``one variable'' projection $P\ci{Q(\bz_j)}$ onto the subspace 
$$
Q(\bz_j) := H^2(\D;E)\cap \Pi(\,\cdot\,,\bz_j) L^2(\D;E)\subset H^2 =H^2(\D; E), 
$$
see Remark \ref{r6.5}.

To prove the second statement of the lemma let us notice that a product of commuting orthogonal projections is an orthogonal projection. Therefore $P=P_{Q_1}P_{Q_2}\ldots P_{Q_n}$ is an orthogonal projection. 

Since for $\xi \in H^2(\D^n;E)\cap \Pi L^2=Q\subset Q_j$
$$
P_{Q_j}\xi =\xi\qquad \forall j=1, 2, \ldots, n, 
$$
we can conclude that 
$$
Q=H^2(\D^n;E)\cap \Pi L^2\subset \ran P. 
$$
On the other hand, since the projections $P_{Q_j}$ commute and $\ran P_{Q_j} = H^2_j\cap \Pi L^2$
$$
\ran P \subset H^2_j\cap \Pi L^2 =Q_j\qquad \forall j=1, 2, \ldots, n,
$$
so 
$$
\ran P \subset \bigcap_{j=1}^n Q_j =\bigcap_{j=1}^n H^2_j\cap \Pi L^2 =H^2\cap \Pi L^2 =Q. 
$$
Therefore $\ran P=Q$, i.e.~$P$ is the orthogonal projection onto $Q$.  
\end{proof}

We can now move onto proving Lemma \ref{l6.2}.

\begin{proof}[Proof of Lemma \ref{l6.2}]
  We will follow the argument  in Lemma \ref{l6.1}. For $\xi\in K$ consider the orthogonal decomposition
$$
\xi= P_{K_1}\xi + P_{Q_1}\xi =: \xi_1 + \xi^1, \qquad \xi_1\in K_1(\D^n;E).
$$
Since $\xi_1\perp\xi^1$, 
$$
\|\xi\|_2^2 =\|\xi_1\|^2_2 +\|\xi^1\|^2_2. 
$$
Decomposing $\xi^1$ as 
$$
\xi^1 = P_{K_2}\xi^1 + P_{Q_2}\xi^1 =: \xi_2 + \xi^2, \qquad \|\xi^1\|_2^2 =\|\xi_2\|^2_2 +\|\xi^2\|^2_2
$$
we get the decomposition of $\xi$
$$
\xi = \xi_1 +\xi_2 + \xi^2, \qquad \xi_j\in K_j,\quad \xi^2= P_{Q_2}P_{Q_1}\xi,  
$$
and
$$
\|\xi\|_2^2 =\|\xi_1\|_2^2 + \|\xi_2\|^2_2 +\|\xi^2\|^2_2 .
$$
Repeating the procedure of decomposing on each step $\xi^k$ using $P_{K_{k+1}}$ we finally obtain
$$
\xi =\xi_1 + \xi_2 +\ldots + \xi_n + \xi^n, \qquad \xi_j\in K_j, \ j=1, 2, \ldots, n, \quad \xi^n = P_{Q_n}\ldots P_{Q_2}P_{Q_1} \xi, 
$$
and 
$$
\|\xi\|_2^2 =\|\xi_1\|_2^2 + \|\xi_2\|^2_2 +\ldots + \|\xi_n\|^2_2 +\|\xi^n\|^2_2 .
$$
But according Lemma \ref{l6.4} $\xi^n=0$, so the lemma is proved. 
\end{proof}

\subsection{Proof of the $H^2$ corona for the polydisk}The idea of the proof is quite simple, we want to reduce everything to one-variable estimates. In the one-variable case we defined the functional $L$ on a dense subset $K$ by
$$
L(\xi) = \int_\D \p  \left[ \La \db\Phi g, \xi \Ra \right] \,d\mu
$$
where $d\mu=\f{2}{\pi}\log\f{1}{\abs{z}}dxdy$, see \eqref{1.3a}. We have also proved (see Proposition \ref{pr1.1}) that the functional $L$ is bounded (in $L^2$ norm on $K$). 

For the polydisk, define (conjugate linear) functionals $L_j$ on $K_j$ by 
$$
L_j (\xi) := \int_{\T^{n-1}} L(\xi(\,\cdot\,,\bz_j)) dm_{n-1}(\bz_j).
$$
Since $\xi(\,\cdot\,,\bz_j) \in K$ for almost all $\bz_j\in \T^{n-1}$ if $\xi\in K_j$ (see Remark \ref{r6.5}) the functionals $L_j$ are well defined and bounded, $\|L_j\|=\|L\|$. Note also, that on a dense set of $\xi$ of the form $\xi=\Pi h$, $h\in (H^2_j)^\perp$ we can represent 
$$
L_j(\xi) = \int_{\T^{n-1}} \int_\D \p_j  \left[ \La \db_j\Phi g, \xi \Ra \right] \,d\mu(z_j) \,dm_{n-1}(\bz_j).
$$

Define a conjugate linear functional $\bL$ on $K$ by decomposing $\xi\in K$ as
\begin{equation}
\label{6.6}
\xi=\xi_1+\xi_2 +\ldots +\xi_n, \qquad \xi_j\in K_j, \quad j=1, 2, \ldots, n
\end{equation}
and putting
$$
\bL(\xi) := \sum_{j=1}^n L_j(\xi_j). 
$$

We will show later that the functional $\bL$ is well defined, i.e.~that it does not depend on the choice of  decomposition of $\xi$ (note that by Lemma \ref{l6.2} one can always find at least one such decomposition).

Assuming for now that $\bL$ is well defined, let us prove Theorem \ref{t0.2} for $p=2$. First of all, by Lemma \ref{l6.2} any function $\xi\in K$ can be decomposed as 
$$
\xi = \sum_{j=1}^n \xi_j, \qquad\text{where } \xi_j\in K_j, \qquad \text{and }\sum_{j=1}^n \|\xi_j\|^2 =\|\xi\|^2.
$$
Therefore, using the fact that $\|L_j\|=\|L\|$ we get for $\xi\in K$
$$
|\bL(\xi)|\le\sum_{j=1}^n \|L_j\|\cdot \|\xi_j\| = \|L\|\sum_{j=1}^n \|\xi_j\| \le \|L\|\sqrt{n} \left(\sum_{k=1}^n \|\xi_j\|^2\right)^{1/2} =\sqrt n\,\|L\|\cdot  \|\xi\|, 
$$
so 
$$
\|\bL\| \le \sqrt n\, \|L\|\le \f{\sqrt{n}C}{\delta^{r+1}}\log\f{1}{\delta^{2r}}
$$
where $C= \sqrt{1+e^2} + \sqrt{e} + \sqrt{2} e \approx 8.38934$ is the constant from Theorem \ref{t0.1}. 

Take $h\in (H^2)^\perp$, and decompose it according to Lemma \ref{l6.1} as
$$
h=\sum_{j=1}^n h_j, \qquad h_j\in (H^2_j)^\perp. 
$$
Denote
$$
\xi:= \Pi h, \qquad \xi_j =\Pi h_j.
$$
Repeating the reasoning with the Green's Formula from the one-variable case we can easily show that
$$
\int_{\T^n} \La \Phi g, h_j\Ra \,dm_{n}(\bz) = \int_{\T^{n-1}} \int_\D \p_j  \left[ \La \db_j\Phi g, \xi_j \Ra \right] \,d\mu(z_j) \,dm_{n-1}(\bz^j)= L_j(\xi_j), 
$$
so
$$
\int_{\T^n} \La \Phi g, h\Ra \,dm_{n}(\bz) = \bL(\Pi h) =\bL(\xi) .
$$

By the Hilbert space version of the Hahn--Banach Theorem the linear functional $\overline\bL$ can be extended to a bounded functional on all of $L^2$, i.e.~we can find $v\in L^2=L^2(\T^n;E)$ such that 
$$
\bL(\xi ) =\int_{\T^n} \La v, \xi\Ra dm_n(\bz)\qquad \forall \xi \in K. 
$$
Replacing $v$ by $\Pi v$ if necessary, one can assume without loss of generality that $v(\bz)\in \ran \Pi(\bz)=\ker F(\bz)$ a.e.~on $\T^n$, so $Fv\equiv 0$ on $\T^n$. 
Since by the construction
$$
\int_{\T^n} \La v, h\Ra \,dm_n(\bz) = \int_{\T^n} \La v, \Pi h\Ra \,dm_n(\bz) = \bL(\Pi h) = \int_{\T^n} \La \Phi g, h\Ra \,dm_n(\bz) \qquad \forall h\in H^2(\D^n;E)^\perp, 
$$
the function $f:=f_0-v := \Phi g -v$ is  analytic. Since $Fv=0$, it satisfies $Ff = Ff_0 =g$, so $f$ is the analytic solution we want to find. \hfill\qed

\subsection{Why the functional $\bL$ is well defined}
Let us consider first the case of the bidisk $\D^2$. To show that $\bL$ is well defined in this case, it is sufficient to show that if 
$$
0=\xi_1 + \xi_2, \qquad \xi_j \in K_j
$$
then $L_1(\xi_1) + L_2(\xi_2)=0$ (take difference of two representations of the same function in $K$). This holds if and only if 
$$
L_1(\xi)=L_2(\xi) \qquad \forall \xi \in K_1\cap K_2. 
$$

Thus, the following lemma shows that $\bL$ is well defined in the case of bidisk $\D^2$. 
\begin{lm}
\label{l6.8}
Let $\xi\in K_1\cap K_2 \subset \Pi L^2(\T^2;E)$. Then 
$$
L_1(\xi)=L_2(\xi)
$$
\end{lm}
\begin{proof}
The proof of this lemma is 
really nothing more than repeated applications of Green's Formula, 
and using that $K_1\cap K_2=\textnormal{clos}(\Pi 
\overline{H^{2}})$ where $\overline{H^{2}}$ are the functions which are anti-holomorphic in both variables.  To see that $K_1\cap K_2=\textnormal{clos}(\overline{\Pi H^{2}})$ we use Lemma \ref{l6.2.1}.  Since $(K_1\cap K_2)^{\perp}=Q_1+Q_2=K_1^{\perp}+K_2^{\perp}=(H^{2}_1+H^{2}_2)\cap\Pi L^{2}$, then by Lemma \ref{l6.2.1} we have the result.

By density we can work with $\xi$ of the form $\xi=\Pi h$ with $h$ anti-holomorphic in 
both variables.Ê So applying Green's Formula twice gives

\begin{eqnarray*}
L_{1}(\xi) & = & \int_{\T}\int_{\D}\p_1\La\pbar_1\Phi g,\xi\Ra d\mu(z_1) 
dm(z_2)\\
 & = & \int_{\T}\int_{\T}\La\Phi g,h\Ra dm(z_1) dm(z_2)\\
 & = & \int_{\T}\int_{\D}\p_2\La\pbar_2\Phi g, \xi\Ra d\mu(z_2) dm(z_1)\\
 & = & L_{2}(\xi).
\end{eqnarray*}
Since this result holds on a dense set of $\xi$, and the functionals $L_1$ and $L_2$ are continuous we have the result for all $\xi\in K_1\cap K_2$.
\end{proof}

For the polydisk the lemma has the following important corollary
\begin{cor}
\label{c6.9}
Let $\xi\in K_j\cap K_k\subset L^2(\T^n;E)$. Then 
$$
L_j(\xi)=L_k(\xi). 
$$
\end{cor}

\begin{proof}
To prove the corollary one needs to apply Lemma \ref{l6.8} to the bidisk in variables $z_j$ and $z_k$ and then integrate the obtained equality over $\T^{n-2}$ (with respect to Lebesgue measure in all other variables). 
\end{proof}

Now we are ready to prove that $\bL$ is well defined. To prove this it is sufficient to show for any representation of $0$ 
\begin{equation}
\label{6.7}
0= \sum_{j=1}^n \xi_j, \qquad \xi_j\in K_j
\end{equation}
the equality 
$$
\sum_{j=1}^n L_j(\xi_j) =0
$$
holds.

We will use induction in $n$. The case $n=2$ is already settled, so let us assume the functional $\bL$ is well defined for the polydisk $\D^{n-1}$. It follows from \eqref{6.7} that 
$$
\xi_n \in K_n \cap (K_1 + K_2 +\ldots + K_{n-1}) = (K_1\cap K_n) + (K_2\cap K_n)+ \ldots + (K_{n-1}\cap K_n), 
$$
so $\xi_n$ can be represented as
$$
\xi_n = \sum_{j=1}^{n-1}\eta_j, \qquad \eta_j\in K_j\cap K_n, \quad j=1, 2, \ldots, n-1. 
$$
On the other hand we know that $\xi_n =-\sum_{j=1}^{n-1}\xi_j$. Using the induction hypothesis and integrating it over $\T$ with respect to $dm(z_n)$ we obtain that
$$
\sum_{j=1}^{n-1} L_j(\eta_j) = - \sum_{j=1}^{n-1} L_j(\xi_j). 
$$
Since $\eta_j\in K_j\cap K_n$, Corollary \ref{c6.9} implies that  $L_j(\eta_j) = L_n(\eta_j)$. Therefore
$$
L_n(\xi_n) = \sum_{j=1}^{n-1} L_n(\eta_j) = \sum_{j=1}^{n-1} L_j(\eta_j) = -\sum_{j=1}^{n-1} L_j(\xi_j), 
$$ 
and so $\sum_{j=1}^{n} L_j(\xi_j) =0$. \hfill \qed

\section{The $H^{p}$ Corona Problem in the Polydisk}
\label{s7}

A simple idea of proving the $H^p$ corona problem in the polydisk is to try to mimic the proof of the $H^2$ case. However, there is a much easier way: just use objects which are already defined, and modify the crucial estimates. 

First of all notice, that replacing the Corona data $F$ and $g$ by $F(r\bz)$ and $g(r\bz)$, $r<1$ and using the standard normal families argument one can assume without loss of generality (as long as we are getting the same uniform estimates on the norm of the solution) that both $F$ and $G$ are holomorphic in a slightly bigger polydisk. So we can always assume that, for example, the right hand side $g$ is not only in $H^p$, but is also bounded, smooth, etc. 

As in the $H^2$ case we first construct a smooth solution $f_0:=\Phi g$, where $\Phi := F^*(FF^*)^{-1}$, of the equation $Ff=g$ and then correct it to be analytic.   To do that it is sufficient to show that the conjugate linear functional $\bL$ introduced in the previous section is $L^q$ bounded, $1/p+1/q=1$, i.e.~that
$$
|\bL(\xi)|\le C\|\xi\|_p
$$
for all $\xi$ of form $\xi= \Pi h$, where $h$ is a trigonometric polynomial in $H^2(\D^n;E)^\perp$. 

If this estimate is proved, the linear functional $\overline \bL$ can be extended by the Hahn--Banach Theorem to a linear functional on $L^q$, so there will exist a function $v\in L^p(\T^n;E)$, $\|v\|_p=\|\bL\|_p$ such that
$$
\bL(\xi) =\int_{\T^n} \La v,\xi\Ra \,dm_n(\bz) \qquad \forall \xi=\Pi h, \quad h\in H^2(\D^n;E)^\perp \cap \operatorname{Pol}.
$$
Again, replacing $v$ by $\Pi v$ we can always assume without loss of generality that $v(\bz) \in \ran \Pi(z) = \ker F(z)$ a.e.~on $\T^n$. As in the previous section, decomposing $h$ as 
$$
h=\sum_{j=1}^n \xi_j, \qquad \xi_j \in H^2_j
$$
($h$ is a trigonometric polynomial, so we can use Lemma \ref{l6.1} here), we can show that $\int_{\T^n} \La \Phi g, h\Ra \,dm_{n}(\bz) = \bL(\Pi h) =\bL(\xi)$ so
$$
\int_{\T^n} \La v, h\Ra \,dm_n(\bz) = \int_{\T^n} \La v, \Pi h\Ra \,dm_n(\bz) = \bL(\Pi h) = \int_{\T^n} \La \Phi g, h\Ra \,dm_n(\bz) \qquad \forall h\in H^2(\D^n;E)^\perp\cap\operatorname{Pol}. 
$$
Therefore, the function $f=f_0-v=\Phi g -v$ is analytic, and it clearly solves the equation $Ff=g$ (on $\T^n$, and therefore on $\D^n$).

\subsection{Main estimates}
Let us introduce some notation. Denote 
$$
K^q:= \clos(\Pi ((H^p)^\perp))\subset \Pi L^q, \qquad Q^q := H^q\cap \Pi L^q, 
$$
so for $K$ and $Q$ introduced in the previous section $K=K^2$ and $Q=Q^2$. Let also 
$$
H^q_j=H^q_j(\D^n;E):= \{f\in L^q (\T^n;E): f(\,\cdot\,,\bz^j) \in H^q(\D;E)\}
$$
be the spaces of functions analytic in variable $z_j$, and let
$$
K^q_j := \clos(\Pi (H^p(\D^n;E)^\perp))\subset \Pi L^q(\T^n;E), \qquad Q^q_j: = H^q(\D^n;E)\cap \Pi L^q(\T^n;E).
$$

To estimate the functional $\bL$ we need the following analogue of Lemma \ref{l6.2} 

\begin{lm}
\label{l7.1-new}
Any function $\xi \in K^q$ can be decomposed as
$$
\xi =\sum_{j=1}^n \xi_j, \qquad \xi_j \in K^q_j, \quad \|\xi_j\|_q\le C(q)^j, 
$$
where $C(q) =1/\sin(\pi/q)$ is the norm of the scalar Riesz Projection $P_+$ from $L^q(\T)$ onto $H^q(\D)$ (note that $C(p) =C(q)$ for $1/p+1/q=1$). 
\end{lm}

Let us show how this lemma implies the estimate for $\bL$. In Section \ref{s5} we have proved the $L^p$ bound for the functional $L$ (in the one-variable case), 
$$
|L(\xi)|\le C(r,\delta) \|\xi\|_p, \qquad C(r,\delta)=C\frac{1}{\delta^{r+1}} \log\frac{1}{\delta^{2r}}, 
$$
where $C= \sqrt{1+e^2} + \sqrt e + \sqrt 2 e$. That would imply the same estimates for the functionals $L_j$ on $L^q(\T^n;E)$, so applying Lemma \ref{l7.1-new} we get
$$
|L(\xi)|\le C(r,\delta) \sum_{j=1}^n \|\xi_j\|_q \le C(r,\delta) \|\xi\|_q \sum_{j=1}^n C(q)^j \le C(r,\delta) n C(q)^n \|\xi\|_q
$$
Recalling that $C(p)=C(q)$ we get the desired estimate of the solution. 

There is a little detail here as the functional $\bL$ was defined initially only on $K^2$.  So formally, if $q<2$ (i.e.~if $p>2$) the functional is not defined on $K^q$. However this is not a big problem and the simplest way of dealing with it is to use the standard approximation arguments. Since the polynomials in $(H^2_j)^\perp\cap \operatorname{Pol}$ are dense in $(H^p_j)^\perp$, the functions of form $\Pi h$, $h\in H^2_j\cap \operatorname{Pol}$ are dense in $K^q_j$. So, approximating functions $\xi_j$ from Lemma \ref{l7.1-new} by functions of this form, we will get the desired estimate. Note, that we are estimating $\bL(\xi)$ on a dense set of functions $\xi=\Pi h$, $h\in (H^2)^\perp \cap \operatorname{Pol}$, so we do not need it be formally defined on $K^q$. 

The main step in proving Lemma \ref{l7.1-new} is the following result that states that in the one-variable case the norm of the orthogonal projections $P_K$ and $P_Q$ in $L^q$ is the same as the norm of the Riesz projection $P_+$ in $L^q$. 

\begin{lm}
\label{l7.2-new}
Let $H^2 = H^2(\D;E)$ and let $K, Q\subset H^2$ be the subspaces defined above in \eqref{6.2a} and \eqref{6.5}. Then for $1<q<\infty$
$$
\|P_K\xi\|_q \le C(q) \|\xi\|_q, \quad \|P_Q\xi \|_q\le C(q) \|\xi\|_q\qquad \forall \xi\in \Pi L^2\cap \Pi L^q
$$
where $C(q)=1/\sin(\pi/q)$ is the norm of the Riesz Projection $P_+$ in $L^q$ (or in $L^p$, $1/p+1/q=1$). 
\end{lm}
Note that since $\Pi L^2\cap\Pi L^q$ is dense in $\Pi L^q$, the projections $P_{K}$ and $P_{Q}$ extend to bounded operators on $\Pi L^q$.
\begin{proof}[Proof of Lemma \ref{l7.2-new}]
Take $\xi\in \Pi L^2\cap \Pi L^q$ and decompose it as 
$$
\xi = P_K \xi + P_Q\xi =: \xi\ci K + \xi\ci Q.
$$
Since $Q$ is a $z$-invariant subspace of $H^2(\D, E)$, by the Beurling--Lax theorem, see \cite{Nikolskii0} it can be represented as 
$Q=\Theta H^2(\D;E_*)$, where $\Theta \in H^\infty(E_*\shto E)$ is an inner function (i.e.~$\Theta(z)$ is an isometry a.e.~on $\T$) and $E_*$ is an auxiliary Hilbert space. So $\xi\ci Q$ can be represented as
$$
\xi\ci Q =\Theta \eta, \qquad \eta\in H^2(E_*)\cap H^q(E_*). 
$$
By duality
$$
\|\xi\ci Q\|_q =\|\eta\|_q = \sup\begin{Sb}h\in L^p\cap L^2:\\ \|h\|_q=1\end{Sb} \left| \int_\T \La \eta , h\Ra \,dm \right| .
$$
Let $h_+= P_+h$. Since $\eta\in H^2$
$$
\int_\T \La \eta , h\Ra \,dm = \int_\T \La \eta , h_+\Ra \,dm = \int_\T \La \Theta \eta , \Theta h_+\Ra \,dm =
\int_\T \La \xi\ci Q  , \Theta h_+\Ra \,dm = \int_\T \La \xi , \Theta h_+\Ra \,dm;
$$
the second equality holds because $\Theta$ is an isometry a.e.~on $\T$, and the last one holds because $\xi_K\in K\perp \Theta h_+$. Therefore, since $\|h_+\|_p\le C(p)\|h\|_p$, we can conclude
$$
\left| \int_\T \La \eta , h\Ra \,dm \right| \le \left| \int_\T \La \xi , \Theta h_+\Ra \,dm\right| \le \|\xi\|_q \|h_+\|_p\le C(p) \|\xi\|_q \|h\|_p
$$
so $\|\xi\ci Q\|_q\le C(p) \|\xi\|_q$. Thus we get the desired estimate for the norm of $P_Q$. 

Since $P_K + P_Q=I$ we can estimate the norm of $P_K$ by $C(p)+1$ for free. Note, that unlike the case of Hilbert spaces, complementary projections in Banach spaces do not necessarily have equal norms. So, to get rid of the $1$ some extra work is needed. 

It is easy to see that $\cap_{n>0}\overline z K =\{0\}$, so the decomposition $\Pi L^2 = K\oplus Q$ implies that the set 
$$
\bigcup_{n>0} \overline z^n Q = \bigcup_{n>0} \overline z^n \Theta  H^2(E_*) 
$$ 
is dense in $\Pi L^2$. Thus $\Pi L^2 = \Theta L^2$, and since $\Theta $ is an isometry a.e.~on $\T$ we can conclude that $K=\Theta ( H^2(E)^\perp) $. Therefore we can represent $\xi\ci K$ as 
$$
\xi\ci K = \Theta \eta, \qquad \eta \in H^2(E_*)^\perp\cap L^q(E_*).
$$
Performing the same calculations as in the case of $\xi\ci Q$, only using $h_-=P_-h$, $P_-=I-P_+$ instead of $h_+$ we get the estimate $\|P_K\|_{L^q}\le \|P_-\|_{L^q}$. But the isometry $\tau$, 
$$
\tau(z^k ) = z^{-k-1}, \qquad k\in \mathbb{Z}
$$
interchanges $H^2$ and $(H^2)^\perp$, and since $\tau$ is an isometry in all $L^p$, we conclude the $\|P_-\|_{L^q}= \|P_+\|_{L^q}$. 
\end{proof}

\begin{cor}
\label{c7.3}
Let $H^{2}=H^{2}(\D^{n};E)$ and let $K_j, Q_j\subset H^{2}$ be the subspaces defined in \eqref{6.3a} and \eqref{6.5}.  Then for $1<q<\infty$ and $1\leq j\leq n$ we have
$$
\|P_{K_j}\xi\|_q \le C(q) \|\xi\|_q, \quad \|P_{Q_j}\xi \|_q\le C(q) \|\xi\|_q\qquad \forall \xi\in \Pi L^2\cap \Pi L^q
$$
where $C(q)=1/\sin(\pi/q)$ is the norm of the (one-dimensional) Riesz Projection $P_+$ in $L^q$ (or in $L^p$, $1/p+1/q=1$). 
\end{cor}

\begin{proof}
This corollary follows directly from Lemma \ref{l7.2-new}.  Since by Remark \ref{r6.5} we can view $P_{K_j}$ and $P_{Q_j}$ as ``one-variable'' operators.  Then we ``freeze'' all variables except the $z_j$ variable and apply Lemma \ref{l7.2-new} and then integrate in the ``frozen'' variables.
\end{proof}
It only remains to prove Lemma \ref{l7.1-new}.
\begin{proof}[Proof of Lemma \ref{l7.1-new}]
Since the projections $P_{K_j}$ extend to bounded operators on $\Pi L^{q}$ we can use Lemma \ref{l6.2} about the decomposition of $\xi\in K^{2}$ and apply Corollary \ref{c7.3} to get the norm of each $\xi_j$.  Indeed, take $\xi\in K^{q}$ and then repeating the proof  of Lemma \ref{l6.2} we can write
$$
\xi=P_{K_1}\xi+P_{Q_1}\xi:=\xi_1+\xi^{1}.
$$

By Corollary \ref{c7.3} we have that $\xi_1\in K_1^{q}$ with $\norm{\xi_1}_{q}\leq C(q)\norm{\xi}_{q}$ and $\norm{\xi^1}_q\leq C(q)\norm{\xi}_q$.  Decomposing $\xi^1$ in the same manner we have
$$
\xi=\xi_1+\xi_2+\xi^2,\quad\xi_j\in K_j^q,\quad \xi^2=P_{Q_2}P_{Q_1}\xi
$$
and by Corollary \ref{c7.3} $\norm{\xi_j}_{q}\leq C(q)^{j}\norm{\xi}_{q}$ because each of $\xi_j$ will have exactly $j$ copies of the one-variable operators $P_K$ or $P_Q$ appearing.  Continuing this decomposition at each step we find
$$
\xi=\xi_1+\xi_2+\ldots+\xi_n+\xi^n,\quad\xi_j\in K_j^q,\quad \xi^n=P_{Q_n}\ldots P_{Q_2}P_{Q_1}\xi,
$$
and $\norm{\xi_j}_q\leq C(q)^j\norm{\xi}_q$ for the same reason as above.  Finally, by Lemma \ref{l6.4} we have that $\xi^n=0$ because $P_{Q_n}\cdots P_{Q_1}=0$ on a dense set.  
\end{proof}

\section{Concluding Remarks}
As was mentioned before, unfortunately this proof does not provide a 
solution to the $H^{\infty}$ corona problem for the polydisk because 
the projections that we used are not bounded.

Additionally, we should point out that a small generalization of the 
problem discussed in this paper could be dealt with using the same techniques
with very little modification necessary to deal with it.  One could 
instead consider the problem where we have
$$
\delta^{2}I\leq FF^{*}\vert_{(\ker F^{*})^{\perp}}\leq I
$$
and $g(z)\in \ran F(z)$ for all $z\in\D$ (or $\D^{n}$ depending on 
the problem one is considering).  This condition on $g$ will ensure 
that $\Phi g$ where $\Phi:=F^{*}(F^{*}F+\epsilon I)^{-1}$ will be well 
defined.  But then the argument can be carried out as above and one 
can then let $\epsilon$ go to zero.  This will be the result of another paper.

\end{document}